\documentclass[12pt]{article}
\usepackage[utf8]{inputenc}
\usepackage{amsmath,amssymb,amsthm}
\renewcommand\le{\leqslant}
\renewcommand\ge{\geqslant}
\newcommand\eps{\varepsilon}
\newcommand\R{\mathbb{R}}

\newtheorem*{lemmaA}{Lemma A}
\newtheorem*{lemmaB}{Lemma B}

\DeclareMathOperator\rank{rank}
\DeclareMathOperator\supp{supp}
\DeclareMathOperator\conv{conv}
\DeclareMathOperator\Span{span}
\DeclareMathOperator\extr{extr}

\newtheorem*{theorem}{Theorem}
\newtheorem{corollary}{Corollary}

\title{On the structure of low-rank matrices that approximate the identity
matrix}
\author{Yuri Malykhin}
\begin{document}
\maketitle

Consider the problem of elementwise approximation of the $N\times N$ identity
matrix by matrices of rank not exceeding given $n$:
\begin{equation}
    \label{approx}
\min_{\rank A\le n}\;\max_{1\le i,j\le N}|A_{i,j}-\delta_{i,j}|.
\end{equation}
The quantity~\eqref{approx} equals to the Kolmogorov width of the octahedron,
$d_n(B_1^N,\ell_\infty^N)$. This width was studied in the
papers~\cite{K74,Gl86,H79}, etc.  Some applications to combinatorics are given
in~\cite{A09}.

The probabilistic method ($A_{i,j}=N^{-1}\langle x^i,x^j\rangle$ with random
$x^1,\ldots,x^N\in\{-1,1\}^n$, see.~\cite{K74}) gives the bound $d_n(B_1^N,\ell_\infty^N)\le
2\sqrt{\ln(N)/n}$; therefore, good approximation of the identity matrix is
possible for $n\asymp\log N$. On the other hand, it is well known
that a nontrivial approximaion requires at least
logarithmic rank. Say, if the error~\eqref{approx} does not exceed $1/3$,
then $N$ columns of the approximation matrix lie inside a ball of radius $5/6$
in an $n$-dimensional subspace of $\ell_\infty^N$; and they are
$1/3$-separated. Standard volume argument shows that $N\le 6^n$.

B.S.~Kashin asked a question: is it true that for matrices of rank
$n\asymp\log N$ that approximate the identity matrix, positive proportion of
elements have absolute value greater than $\gamma$, with $\gamma\asymp n^{-1/2}$?

Given a matrix $A\in\R^{N\times N}$, we define the distribution function
$F^*_A(\gamma) := N^{-2}|\{(i,j)\colon |A_{i,j}|>\gamma\}|$.

\begin{theorem}
    Let $A$ be a $N\times N$ matrix of rank $n$, such that
    $$
    \max_{1\le i,j\le N}|A_{i,j}-\delta_{i,j}|\le 1/3.
    $$
    Then the inequality holds
    $$
    F_A^*(\gamma) \ge c\frac{\log N}{n\log(2+\frac{n}{\log N})},
    \quad\mbox{where $\gamma := c \max\{n^{-3/2}\log N,\,n^{-1}\}$}.
    $$
    Here $c$ is an absolute positive constant.
\end{theorem}

\begin{corollary}
    In the conditions of the Theorem, if $n\le K\log N$ for some $K$, then at least
    $c(K)N^2$ elements of the matrix $A$ have absolute value $>c(K)n^{-1/2}$,
    where $c(K)>0$ and depends only on $K$.
\end{corollary}
This corollary gives an affirmative answer to Kashin's question.
The probabilistic method shows that the order $n^{-1/2}$ is sharp.

\begin{corollary}
    In the conditions of the Theorem, the proportion of nonzero elements of $A$ is at
    least $c \log(N)/(n\log(2+n/\log N))$.
\end{corollary}
It is easy to construct a matrix that approximates the identity but the
proportion of nonzero elements is $\asymp \log(N)/n$.

We will use the following statements.
\begin{lemmaA}[\cite{G95}]
    Let $K\subset\R^n$ be a symmetric convex body, let $D\supset K$ be its minimal
    volume ellipsoid, $|\cdot|_D$ is the Euclidean norm corresponding to $D$.
    Then for any $\eps\in(0,1)$ there exist vectors
    $x^1,\ldots,x^k\in K\cap\partial D$, $k\ge n(1-\eps)$, such that for any
    reals $t_1,\ldots,t_k$ the inequality holds
    \begin{equation}
        \label{ineq_ellispoid}
    |\sum_{m=1}^k t_mx^m|_D \ge c_0 \eps n^{-1/2} \sum_{m=1}^k |t_m|.
    \end{equation}
\end{lemmaA}
This lemma was proven in~\cite{ST89} with the coefficient $\eps^{3/2}n^{-1/2}$;
in~\cite{G95} the result was improved to $\eps n^{-1/2}$.
Lemma was used to bound the Banach--Mazur distance from an arbitrary 
$n$-dimensional space to $\ell_\infty^n$. Giannopoulos~\cite{G95} proved that
$d(X_n,\ell_\infty^n)\le Cn^{5/6}$.

\begin{lemmaB}
    Let $K\subset\R^n$ be a symmetric convex body. There exist
    $x^1,\ldots,x^n\in\extr K$, such that any
    $x\in K$ is represented as $x=\sum_{m=1}^n t_mx^m$ with some
    $t_1,\ldots,t_n\in[-1,1]$.
\end{lemmaB}
The basis $x^1,\ldots,x^n$ is just an Auerbach basis. One can construct it by a
volume maximization~\cite[Ch.1]{LT} so we can take $x^m\in\extr K$.

Let us prove the Theorem.
\begin{proof}
    Let us say that elements $|A_{i,j}|>\gamma$ are ``large'', and let us call
    the proportion of large elements by ``density''.

The matrix $A$ has density $F_A^*(\gamma)$. Let $\varkappa$ be the maximal
    density of columns of $A$. We can assume that
    $\varkappa=O(F_A^*(\gamma))$ and prove a lower bound for $\varkappa$.
    Indeed, there are at most half columns with density
$>2F_A^*(\gamma)$; one can delete columns and rows with that indexes and pass to
a half-sized submatrix.

Let $v^1,\ldots,v^N$ be the columns of $A$; $V := \Span\{v^j\}_{j=1}^N$~---
    the $n$-dimensional subspace in $\R^N$; $K := \conv\{\pm v^j\}_{j=1}^N$;
    $D\subset V$ is the minimal volume ellipsoid of $K$;
    $\|\cdot\|_K$~--- the norm in $V$ with unit ball $K$;
    $\langle\cdot,\cdot\rangle_D$ and $|\cdot|_D$~--- scalar product and norm
    in $V$ that correspond to the ellipsoid $D$.

    Fix $\eps\in(0,1)$; we will specify its value later. By Lemma A there exist
    vectors $x^1,\ldots,x^k\in K\cap\partial D$, $k\ge n(1-\eps)$, such that the
    inequality~\eqref{ineq_ellispoid} holds.
    We have $K\cap\partial D \subset \extr K \subset \{\pm
    v^j\}_{j=1}^N$, so w.l.o.g.
    $x^1=v^{j_1},\ldots,x^k=v^{j_k}$ with some indexes
    $j_1,\ldots,j_k$. (Soon we will use this.)

    Let us add to vectors $\{x^m\}_1^k$ some vectors $y^1,\ldots,y^{n-k}\in V$
    orthogonal to $\{x^m\}_1^k$ (with respect to $\langle\cdot,\cdot\rangle_D$)
    to make a basis of $V$. Now we can expand columns of our matrix in
    this basis:
$$
v^j = \sum_{m=1}^k t^j_m x^m + \sum_{l=1}^{n-k}s^j_l y^l,\quad j=1,\ldots,N. 
$$
We have $v^j\in K$, so
$$
1 \ge \|v^j\|_K \ge \|v^j\|_D \ge |\sum_{m=1}^k t^j_m x^m|_D
    \ge c_0\eps n^{-1/2} \sum_{m=1}^k |t^j_m|.
$$
Therefore
$$
    A_{i,j} = (v^j)_i = \sum_{m=1}^k t^j_m x^m_i + \sum_{l=1}^{n-k}s^j_l y^l_i
    = \langle x_i, t^j \rangle + \langle y_i, s^j \rangle,
    \quad \|t^j\|_1 \le (c_0\eps)^{-1}n^{1/2}.
$$
where $x_i,t^j\in\R^k$; $y_i,s^j\in\R^{n-k}$.

    Vectors $x_i=(x^1_i,\ldots,x^k_i)=(v^{j_1}_i,\ldots,v^{j_k}_i)$ are rows of
    the matrix that consists of columns of $A$ with indexes
$j_1,\ldots,j_k$. This matrix has density $\le\varkappa$.
Hence, at least half of its rows has density
$\le 2\varkappa$; let us denote by $I$ the set of indexes of these rows.
So, the density of $x_i$ is at most $2\varkappa$ for $i\in I$.

Write $x_i$ as $x_i = w_i + z_i$, where $w_i$ contains only large coordinates of
$x_i$ and $z_i$ contains all other coordinates. We have
$$
|\langle x_i,t^j\rangle - \langle w_i,t^j\rangle| \le \|t^j\|_1 \cdot \gamma \le
    (c_0\eps)^{-1}n^{1/2}  \gamma \le 1/15
$$
provided that $\gamma/\eps \le 1/(15c_0n^{1/2})$.
Consider the matrix
$$
B_{i,j} = \langle w_i, t^j \rangle + \langle y_i, s^j \rangle,
\quad i,j\in I.
$$
We will use three properties:
the number of nonzero elements in each vector
$w_i$ is at most $m :=
\lfloor 2\varkappa k\rfloor$; the dimension of $y_i,s^j$ is at most
$n\eps$; the matrix $B$ approximates the identity:
$|B_{i,j}-\delta_{i,j}|\le 2/5$.

Consider the rows of $B$. They are $1/5$-separated in
$\ell_\infty^I$. On the other hand, we can construct an $(1/11)$--net in
$\ell_\infty^I$ for the set of rows in the following way. For any $i\in I$
the set $\Lambda = \supp w_i$ has at most $m$ elements, so there are at most
$(en/m)^m$ variants for $\Lambda$. (If $m=0$,
we put $(en/m)^m=1$.)

Fix a set $\Lambda$ and consider the submatrix $B_\Lambda$ of $B$ that consists
of rows $i$ such that $\supp w_i=\Lambda$. We have $B_{i,j}=\langle
\left.w_i\right|_\Lambda, \left.t^j\right|_\Lambda\rangle + \langle
y_i,s^j\rangle$ for this submatrix; hence $\rank B_\Lambda\le m+n\eps$. Rows of
$B_\Lambda$ lie in a ball (of fixed diameter) in the space of dimension at most
$m+n\eps$; the size of the net for such a ball does not exceed $\exp(C(m+n\eps))$.

The size of the net is at least the number of rows in $B$, so
$$
(en/m)^m \cdot \exp(C(m+n\eps)) \ge |I| \ge N/2.
$$
Take the logarithm: $m\ln(en/m) + Cm + Cn\eps \ge \ln(N/2)$.
Define
$$
\eps:=\min\left\{\frac12,\frac{\ln(N/2)}{2Cn}\right\},
$$
then
$m\ln(C_1 n/m) \ge \ln(N/2)/2$. It follows that $m\ne0$.
Moreover, $k\ge n(1-\eps)\ge n/2$; $2\varkappa n \ge m=\lfloor 2\varkappa k\rfloor \ge
\varkappa k \ge \varkappa n/2$, therefore:
$$
\varkappa \ln(2C_1/\varkappa) \ge \frac{\ln(N/2)}{4n}.
$$
This proves the required bound for $\varkappa$ with $\gamma \asymp n^{-1/2}\eps$. We know that
$A$ approximates the identity, so $n\ge c_3\log N$ and hence
$\gamma\asymp n^{-3/2}\log N$ in both cases ($\eps=1/2$, $\eps<1/2$).

For $\gamma \asymp n^{-1}$ one should use Lemma B; everything is analogous.
\end{proof}

Let us mention another approach to the problem of estimating the number of large
elements of a matrix for given parameters 
$N,n,\gamma$. Consider the minimal $M$ such that
$d_n(B_1^M,\ell_\infty^M)>\gamma$.
Let $A$ be a matrix of rank $n$ that approximates the $N\times N$
identity matrix. W.l.o.g. $A$ is symmetric. Consider the graph $G$ with $N$
vertices and edges between pairs $i,j$ such that
$|A_{i,j}|\le\gamma$. Then $G$ does not contain complete subgraphs with $M$
vertices (otherwise an approximation of $M\times M$ identity matrix with an
error $\le\gamma$ occurs).
The number of edges in $G$ is bounded via Turan's theorem by
$(1-(M-1)^{-1})N^2/2$.

\end{document}